\documentclass[10pt]{amsart}

\usepackage{amsmath}
\usepackage{amssymb}
\usepackage{bm}
\usepackage{graphicx}
\usepackage{psfrag}
\usepackage{color}
\usepackage{hyperref}
\hypersetup{colorlinks=true, linkcolor=blue, citecolor=magenta, urlcolor=wine}
\usepackage{url}
\usepackage{algpseudocode}
\usepackage{fancyhdr}
\usepackage{mathtools}
\usepackage{tikz-cd}
\usepackage{xy}
\input xy
\xyoption{all}
\usepackage{stmaryrd}
\usepackage{calrsfs}
\usepackage{enumitem}

\newcommand{\ZZ}{\mathbb{Z}}

\newcommand{\FF}{\mathbb{F}}

\newcommand{\vp}{\varphi}
\newcommand{\id}{\mathrm{id}}

\newcommand{\Span}{\operatorname{Span}}
\newcommand{\ten}{\otimes}

\newcommand{\End}{\operatorname{End}}
\newcommand{\al}{\alpha}
\newcommand{\sub}{\subset}

\newcommand{\abs}[1]{{\left\lvert #1\right\rvert}}

\newcommand{\br}[1]{{\left\{ #1 \right\}}}
\newcommand{\pr}[1]{{\left( #1 \right)}}

\newcommand{\ob}[1]{\mkern 1.5mu\overline{\mkern-1.5mu#1\mkern-1.5mu}\mkern 1.5mu}

\newcommand{\db}[1]{\llbracket #1 \rrbracket}

\newcommand{\rst}[1]{\!\mid_{#1}}
\newcommand{\up}[1]{^{\pr{#1}}}

\newenvironment{enum}
    {
    \begin{enumerate} [label = \bf{(\alph*)}, topsep = 0pt, itemsep = -1ex, partopsep = 1ex, parsep = 1ex]
    }
    {
    \end{enumerate}
    \vspace{1mm}
    }

\newcommand{\INDEC}{\operatorname{INDEC}}
\newcommand{\noalpha}{\INDEC[X_1, \ldots, X_M]}

\newtheorem*{maintheorem*}{Main Theorem}

\newtheorem*{theorem*}{Main Theorem}

\theoremstyle{definition}

\numberwithin{equation}{section}

\keywords{Quiver, Shur representation, 3-cnf formula, NP-complete problem}

\subjclass[2020]{16G20, 68Q25}

\begin{document}
	
\mbox{}
\title{The Quiver Problem is NP Complete}

\author{Victor Kac}
\address{Department of Math \\ MIT \\ Cambridge, Massachusetts 02139}
\email{kac@math.mit.edu}

\author{Bangzheng Li}
\address{Department of Math \\ MIT \\ Cambridge, Massachusetts 02139}
\email{liben@mit.edu}

\date{\today}
	
\begin{abstract}
    We prove that the quiver problem is NP complete.
\end{abstract}

\bigskip
\maketitle

%%%%%%%%%%%
%%%%%%%%%%%
\section{Introduction}
\label{sec: intro}

Recall that a \underline{representation} $R$ of a quiver (i.e., oriented graph) $Q$ over a field $\FF$ is a collection of finite-dimensional vector spaces $U_v$ over $\FF$, labeled by vertices of $Q$, and linear maps $R_{v, w} : U_v \to U_w$, labeled by oriented edges $v \to w$ of $Q$.
The representation $R$ is called absolutely indecomposable if it does not decompose in a direct sum of two non-zero representations over the algebraic closure of $\FF$.
Choosing bases of the vector spaces $U_v$, $R$ is given by a collection of matrices $R_{v, w}$ over $\FF$ of size $n_v \times n_w$ for each oriented edge $v \to w$, where $n_v = \dim U_v$.
Given such a collection of $\sum_{v \to w} u_v u_w$ elements of $\FF$, the output is YES if $R$ is absolutely indecomposable, and NO otherwise.
Call this problem INDEC.
This is a P-problem if $\FF$ is a finite field, according to \cite[Theorem 2]{K19}.

Define a generalization of INDEC for a non-negative integer $M$, denoted $\INDEC[X_1, \ldots, X_M]$, where the entries of matrices $R_{v, w}$ are either elements of $\FF$ or one of the indeterminats $X_1, \ldots, X_M$.
Say YES for this problem if there exist $x_1, \ldots, x_M \in \FF$ that we can substitute $X_1, \ldots, X_M$, such that the resulting INDEC problem is YES, and NO otherwise.
We call this the \underline{quiver problem}.
It follows from \cite[Theorem 2]{K19} that this is an NP problem if $\FF$ is a finite field.

In the present paper we prove that the problem $\INDEC[X_1, \ldots, X_M]$ is NP complete if $\FF = \FF_2$.

In more detail, let $\Theta = E_1 \wedge \cdots \wedge E_L$ be a 3-cnf Boolean formula \cite[Chapter 7]{S12}, where $E_k = \beta_1\up{k} \vee \beta_2\up{k} \vee \beta_3\up{k}, k = 1, \ldots, L$, are clauses, $\beta_j\up{k} \in \br{X_1, \ldots, X_M, \ob{X_1}, \ldots, \ob{X_M}}$, and the indeterminates $X_i$ take values 0 or 1 in $\FF_2$, so that $\ob{X_i}$ takes values 1 or 0 respectively.

We will construct a quiver $Q_\Theta$ and its representation $R_{\Theta}(X_1, \ldots, X_M)$ from $\noalpha$, such that $(x_1, \ldots, x_M) \in \FF_2^M$ makes $\Theta$ evaluate to 1 if and only if the representation $R_\Theta(x_1, \ldots, x_M)$ is absolutely indecomposable.
This will imply that the problem $\noalpha$ over $\FF_2$ is NP complete, since the problem 3-SAT is NP complete (see \cite[Chapter 7]{S12}).

Recall that an \underline{endomorphism} of a representation $R$ is a collection of endomorphisms $a_v : U_v \to U_v$, such that $R_{v, w} a_v = a_w R_{v, w}$.
For example each $c \in \FF$ defines an obvious endomorphism $a_v = c I_{U_v}$ of $R$.
Recall that $R$ is called a \underline{Schur representation} if its only endomorphisms are the obvious ones.
We show that, moreover, $(x_1, \ldots, x_M) \in \FF_2^M$ makes $\Theta$ evaluate to 1, then $R_\Theta(x_1, \ldots, x_M)$ is a Schur representation, and if it makes $\Theta$ evaluate to 0, then $R_\Theta(x_1, \ldots, x_M)$ decomposes in a direct sum of two non-zero representations.

\bigskip

\section{Construction of the quiver $Q_\Theta$}
\label{sec: construct quiver}

Here is an example of quiver corresponding to the formula $(X_1 \vee X_2 \vee \ob{X_3}) \wedge (X_2 \vee X_3 \vee X_5) \wedge (X_3 \vee \ob{X_4} \vee \ob{X_5})$ (so $M = 5, L = 3$), where the clauses $X_1 \vee X_2 \vee \ob{X_3}, X_2 \vee X_3 \vee X_5$, and $X_3 \vee \ob{X_4} \vee \ob{X_5}$ correspond to vertices $w_1, w_2$, and $w_3$ respectively.

% https://q.uiver.app/#q=WzAsMTMsWzIsMCwid18yIl0sWzEsMSwid18xIl0sWzMsMSwid18zIl0sWzAsNCwidl8xIl0sWzEsMywidl8yIl0sWzMsMywidl8zIl0sWzQsNCwidl80Il0sWzIsNSwidl81Il0sWzAsOCwidV8xIl0sWzEsNywidV8yIl0sWzMsNywidV8zIl0sWzQsOCwidV80Il0sWzIsOSwidV81Il0sWzgsOV0sWzksMTBdLFsxMCwxMV0sWzExLDEyXSxbMTIsOF0sWzgsM10sWzksNF0sWzEwLDVdLFsxMSw2XSxbMywxLCIiLDAseyJvZmZzZXQiOi0xfV0sWzEsM10sWzQsMV0sWzEsNCwiIiwwLHsib2Zmc2V0IjotMX1dLFs1LDFdLFsxLDUsIiIsMCx7Im9mZnNldCI6LTF9XSxbNCwwXSxbMCw0LCIiLDAseyJvZmZzZXQiOi0xfV0sWzUsMF0sWzAsNSwiIiwwLHsib2Zmc2V0IjotMX1dLFs3LDBdLFswLDcsIiIsMCx7Im9mZnNldCI6LTF9XSxbNSwyXSxbMiw1LCIiLDAseyJvZmZzZXQiOi0xfV0sWzYsMl0sWzIsNiwiIiwwLHsib2Zmc2V0IjotMX1dLFs3LDJdLFsyLDcsIiIsMCx7Im9mZnNldCI6LTF9XSxbMTIsN11d
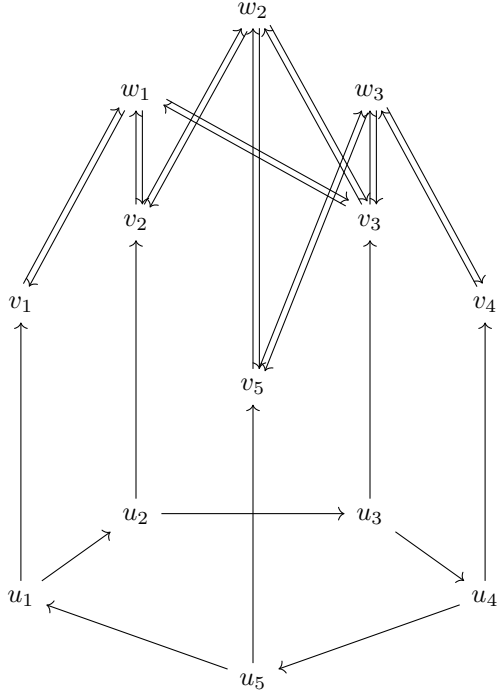
\[\begin{tikzcd}
	&& {w_2} \\
	& {w_1} && {w_3} \\
	\\
	& {v_2} && {v_3} \\
	{v_1} &&&& {v_4} \\
	&& {v_5} \\
	\\
	& {u_2} && {u_3} \\
	{u_1} &&&& {u_4} \\
	&& {u_5}
	\arrow[shift left, from=1-3, to=4-2]
	\arrow[shift left, from=1-3, to=4-4]
	\arrow[shift left, from=1-3, to=6-3]
	\arrow[shift left, from=2-2, to=4-2]
	\arrow[shift left, from=2-2, to=4-4]
	\arrow[from=2-2, to=5-1]
	\arrow[shift left, from=2-4, to=4-4]
	\arrow[shift left, from=2-4, to=5-5]
	\arrow[shift left, from=2-4, to=6-3]
	\arrow[from=4-2, to=1-3]
	\arrow[from=4-2, to=2-2]
	\arrow[from=4-4, to=1-3]
	\arrow[from=4-4, to=2-2]
	\arrow[from=4-4, to=2-4]
	\arrow[shift left, from=5-1, to=2-2]
	\arrow[from=5-5, to=2-4]
	\arrow[from=6-3, to=1-3]
	\arrow[from=6-3, to=2-4]
	\arrow[from=8-2, to=4-2]
	\arrow[from=8-2, to=8-4]
	\arrow[from=8-4, to=4-4]
	\arrow[from=8-4, to=9-5]
	\arrow[from=9-1, to=5-1]
	\arrow[from=9-1, to=8-2]
	\arrow[from=9-5, to=5-5]
	\arrow[from=9-5, to=10-3]
	\arrow[from=10-3, to=6-3]
	\arrow[from=10-3, to=9-1]
\end{tikzcd}\]

% \begin{center}
%     \includegraphics[width = 0.8\textwidth]{Pics/QProbP1.jpg}
% \end{center}

We will first do the following preprocessing to make all variables in a clause $E_k$ to be distinct, i.e., there is no two $X_i$'s or a $X_i$ and a $\ob{X_i}$ appearing at the same time in $E_k$.
First, if there is a $X_i$ and a $\ob{X_i}$ appearing at the same time, then this clause always evaluate to true, so we can simply delete it.

Now we discuss the case where (at least) two $X_i$ appears at the same time in $E_k$.
First note that we can force a variable $z$ to be true by including clauses $(z \vee b \vee c) \wedge (z \vee b \vee \ob c) \wedge (z \vee \ob b \vee c) \wedge (z \vee \ob b \vee \ob c)$, where $b$ and $c$ are new variables.
Then for a clause $X_i \vee X_i \vee Y$, where $X_i \ne Y$, we can simply modify it to $X_i \vee z \vee Y$ ($z$ is a new variable) then force $z$ to be true.
For a clause $X_i \vee X_i \vee X_i$, we can simply replace it by the clauses that forces $X_i$ to be true.
The same argument works for two $\ob{X_i}$'s appearing in the same clause.

With the assumption that all variables in a clause are distinct, we describe our construction of the quiver.
There are $2M + L$ vertices, labeled by $u_1, \ldots, u_M, v_1, \ldots, v_M, w_1, \ldots, w_L$.
To give the arrows between them, for each $i = 1, 2, \ldots, M$, let
\[
\Omega_i^+ = \br{k \in \db{1, L} : X_i \text{ appears in } E_k} = \br{k \in \db{1, L} : \beta_j\up{k} = X_i \text{ for some } j \in \db{1, 3}} \sub \db{1, L},
\]
where $\db{a, b} = \br{k \in \ZZ : a \le k \le b}$.
Similarly, define
\[
\Omega_i^- = \br{k \in \db{1, L} : \ob{X_i} \text{ appears in } E_k} = \br{k \in \db{1, L} : \beta_j\up{k} = \ob{X_i} \text{ for some } j \in \db{1, 3}} \sub \db{1, L}
\]
and
\[
\Omega_i = \Omega_i^+ \cup \Omega_i^- \sub \db{1, L}.
\]

For each $i \in \db{1, M}$, there is an arrow from $u_i$ to $u_{i + 1}$, where we identify $u_{M + 1} = u_1$.

For each $i \in \db{1, M}$, there is an arrow from $u_i$ to $v_i$.

For each $i \in \db{1, M}$ and $k \in \Omega_i$, there is an arrow from $v_i$ to $w_k$ and an arrow from $w_k$ to $v_i$.

There is no other arrows and we have constructed our quiver $Q_\Theta$.

\bigskip

\section{Construction of the representation $R_\Theta(X_1, \ldots, X_M)$}
\label{sec: construct rep}

For each $u_i (i = 1, \ldots, M)$, we assign a 1-dim vector space $U_i = \FF$.

For each $v_i (i = 1, \ldots, M)$, we assign a $n_i := \pr{\abs{\Omega_i} + 1}$-dim vector space $V_i = \FF^{n_i}$.

For each $w_k (k = 1, \ldots, L)$, we assign a 1-dim vector space $W_k = \FF$.

We now describe the linear map for each arrow between them.

\subsection{Bottom horizontal arrows}

For each $i \in \db{1, M}$, take a linear map from $U_i = \FF$ to $U_{i + 1} = \FF$, which is the identity map $\id : \FF \to \FF$.
Here we assume $U_{M + 1} = U_1$.

\subsection{Arrows from $U_i$ to $V_i$}

For each $i \in \db{1, M}$, take a linear map from $U_i = \FF$ to $V_i = \FF^{n_i}$ sending $1 \in \FF$ to $(1, X_i, 0, 0, \ldots, 0) \in \FF[X_1, \ldots, X_M] \ten_\FF V_i$.

We denote by $T_i$ (meaning True) the vector $(1, 1, 0, 0, \ldots, 0) \in V_i$, and $F_i$ (meaning False) the vector $(1, 0, 0, \ldots, 0) \in V_i$.

For each $k \in \Omega_i$, we select a vector $e_i(k) \in V_i$ such that the set (of cardinality $n_i + 1$)
\[
\br{e_i(k) : k \in \Omega_i} \cup \br{T_i, F_i}
\]
has the property that all its subsets of cardinality $n_i$ form a basis in $V_i$.
This can be done by first assigning all but one of the $e_i(k)$'s together with $\br{T_i, F_i}$ to form a basis for $V_i$, then take the last $e_i(k)$ to be the sum of vectors from this basis.

\subsection{$+$ type arrows from $V_i$ to $W_k$}

For each $i \in \db{1, M}$ such that $k \in \Omega_i^+$, take a linear map from $V_i = \FF^{\abs{\Omega_i} + 1}$ to $W_k = \FF$ sending the set 
\[
\br{F_i} \cup \br{e_i(k') : k' \in \Omega_i - \br{k}} \sub V_i
\]
to $0 \in \FF$ and $T_i \in V_i$ to $1 \in \FF$.
This uniquely determines the map (in polynomial time) from the choice of $e_i(k)$'s.

\subsection{$-$ type arrows from $V_i$ to $W_k$}

For each $i \in \db{1, M}$ such that $k \in \Omega_i^-$, take a linear map from $V_i = \FF^{\abs{\Omega_i} + 1}$ to $W_k = \FF$ sending the set 
\[
\br{T_i} \cup \br{e_i(k') : k' \in \Omega_i - \br{k}} \sub V_i
\]
to $0 \in \FF$ and $F_i \in V_i$ to $1 \in \FF$.
This uniquely determines the map (in polynomial time) from the choice of $e_i(k)$'s.

\subsection{Arrows from $W_k$ to $V_i$}

For each $i \in \db{1, M}$ such that $k \in \Omega_i$, take a linear map from $W_k = \FF$ to $V_i = \FF^{n_i}$ sending $1 \in \FF$ to $e_i(k) \in V_i$.

\bigskip

\section{If $(x_1, \ldots, x_M) \in \FF^M$ makes $\Theta$ evaluate to 0 (False), then $R_\Theta(x_1, \ldots, x_M)$ is not indecomposable}
\label{sec: eval 0 is not indecomp}

Then there exists $k \in \db{1, L}$ such that $E_k$ evaluate to 0 (False) under the current selection of $x_1, \ldots, x_M$.

Then we claim $R(x_1, \ldots, x_M) = R_1 \oplus R_2$, where
\[
R_1 = W_k \oplus \bigoplus_{\substack{i \in \db{1, M} \\ k \in \Omega_i}} \FF e_i(k) \sub W_k \oplus \bigoplus_{\substack{i \in \db{1, M} \\ k \in \Omega_i}} V_i.
\]
and
\[
R_2 = \bigoplus_{i = 1}^M U_i \oplus \bigoplus_{i = 1}^M V_i' \oplus \bigoplus_{k' \ne k} W_{k'} \sub \bigoplus_{i = 1}^M U_i \oplus \bigoplus_{i = 1}^M V_i \oplus \bigoplus_{k' \ne k} W_{k'},
\]
where
\[
V_i' = \begin{cases}
    V_i & k \not\in \Omega_i \\
    \FF T_i \oplus \bigoplus_{k' \in \Omega_i - \br{k}} \FF e_i(k') & k \in \Omega_i \text{ and } x_i = 1 \\
    \FF F_i \oplus \bigoplus_{k' \in \Omega_i - \br{k}} \FF e_i(k') & k \in \Omega_i \text{ and } x_i = 0
\end{cases}.
\]

First, it is easy to see that $R(x_1, \ldots, x_M) = R_1 \oplus R_2$ as vector spaces.
Now we check each $R_1$ and $R_2$ are preserved under the arrows.
\begin{enum}
    \item The bottom horizontal arrows between $U_i$'s are obviously preserved.
    \item The vertical arrows from $U_i$ to $V_i$ are preserved from the selection of $V_i'$.
    \item The vertical arrows from $W_k$ to $V_i$'s are preserved because for each possible codomain $V_i (k \in \Omega_i)$ we have $e_i(k)$ is in $R_1$.
    \item The vertical arrows from $W_{k'}$ ($k' \ne k$) to $V_i$'s are preserved because $V_i$ is in $R_2$ if $k \not\in \Omega_i$ and $e_i(k')$ is in $R_2$ if $k \in \Omega_i$.
    \item The vertical arrows from $V_i$ to $W_k$ (hence $k \in \Omega_i$) are preserved because if we write $f : V_i \to W_k$, then $f(e_i(k')) = 0$ for all $k' \ne 0$.
    Also, when $x_i = 1$, then because $E_k$ is evaluated to $0$ we see $k \in \Omega_i^-$, hence $f$ maps $T_i \in V_i$ to $0 \in W_k$.
    When $x_i = 0$, then because $E_k$ is evaluated to $0$ we see $k \in \Omega_i^+$, hence $f$ maps $F_i \in V_i$ to $0 \in W_k$.
    \item The vertical arrows from $V_i$ to $W_{k'} (k' \ne k)$ are preserved because each $e_i(k) \in V_i$ is mapped to $0$ in $W_{k'}$.
\end{enum}

This shows that $R_\Theta(x_1, \ldots, x_M)$ is not indecomposable, as desired.

\bigskip

\section{$R_\Theta(x_1, \ldots, x_M)$ is a Schur representation if $(x_1, \ldots, x_M)$ makes $\Theta$ evaluate to 1}
\label{sec: eval 1 is schur}

Let $\Theta$ be a satisfiable 3-cnf formula.
We will now show that the root we constructed above is a Schur root.
More specifically, we will show that if $(x_1, \ldots, x_M) \in \FF^M$ makes $\Theta$ evaluate to 1, then the representation $R = R_\Theta(x_1, \ldots, x_M)$ satisfies $\End R = \br{c \cdot \id : c \in \FF}$.
In particular, $R$ is absolutely indecomposable.

To see this, pick any $\vp \in \End R$.
Since $\br{c \cdot \id} \sub \End R$, after subtracting suitable multiple of $\id$ we can assume $\vp\rst{U_1} = 0$.
Then we claim $\vp\rst{U_2} = 0$ since we have commutative diagram
% https://q.uiver.app/#q=WzAsNCxbMCwwLCJVXzEiXSxbMSwwLCJVXzEiXSxbMCwxLCJVXzIiXSxbMSwxLCJVXzIiXSxbMCwxLCJcXHZwXFxyc3R7VV8xfSA9IDAiXSxbMSwzLCJcXHNpbSJdLFswLDIsIlxcc2ltIiwyXSxbMiwzLCJcXHZwXFxyc3R7VV8yfSIsMl1d
\[\begin{tikzcd} [column sep = large]
	{U_1} & {U_1} \\
	{U_2} & {U_2}
	\arrow["{\vp\,\rst{U_1} = 0}", from=1-1, to=1-2]
	\arrow["\sim"', from=1-1, to=2-1]
	\arrow["\sim", from=1-2, to=2-2]
	\arrow["{\vp\,\rst{U_2}}"', from=2-1, to=2-2]
\end{tikzcd}\]
Repeating such argument we obtain $\vp\rst{U_i} = 0$ for all $i$.

Now because of the arrows from $U_i$ to $V_i$, for each $i \in \db{1, M}$ such that $x_i = 1$, we have $\vp\rst{V_i}(T_i) = 0$.
This is because we have commutative diagram
% https://q.uiver.app/#q=WzAsNCxbMCwwLCJVX2kiXSxbMSwwLCJVX2kiXSxbMCwxLCJWX2kiXSxbMSwxLCJWX2kiXSxbMCwxLCJcXHZwfF97VV9pfSA9IDAiXSxbMiwzLCJcXHZwfF97Vl9pfSJdLFswLDJdLFsxLDNdXQ==
\[\begin{tikzcd} [column sep = large]
	{U_i} & {U_i} \\
	{V_i} & {V_i}
	\arrow["{\vp|_{U_i} = 0}", from=1-1, to=1-2]
	\arrow[from=1-1, to=2-1]
	\arrow[from=1-2, to=2-2]
	\arrow["{\vp|_{V_i}}", from=2-1, to=2-2]
\end{tikzcd}\]
if we start with $1 \in U_i$ from the upper left, then since $1 \in U_i$ is mapped to $T_i$ under map $U_i \to V_i$, we see $\vp\rst{V_i}(T_i) = 0$.
Similarly, for each $i \in \db{1, M}$ such that $x_i = 0$, we have $\vp\rst{V_i}(F_i) = 0$.

Now for each $k \in \db{1, L}$, there must be some $X_i$ appearing in $E_k$ with $x_i = 1$ or some $\ob{X_i}$ appearing in $E_k$ with $x_i = 0$ (because $(x_1, \ldots, x_M) \in \FF^M$ makes $\Theta$ evaluate to 1, so makes each $E_k$ evaluate to 1).
In the former case, $\vp\rst{V_i}(T_i) = 0$, so because of the arrow from $V_i$ to $W_k$ we see $\vp\rst{W_k} = 0$.
This is because we have commutative diagram
% https://q.uiver.app/#q=WzAsNCxbMCwwLCJXX2siXSxbMSwxLCJWX2kiXSxbMSwwLCJXX2siXSxbMCwxLCJWX2kiXSxbMCwyLCJcXHZwfF97V19rfSJdLFszLDEsIlxcdnB8X3tXX2t9Il0sWzAsM10sWzIsMV1d
\[\begin{tikzcd}
	{V_i} & {V_i} \\
	{W_k} & {W_k}
	\arrow["{\vp|_{V_i}}", from=1-1, to=1-2]
	\arrow[from=1-1, to=2-1]
	\arrow[from=1-2, to=2-2]
	\arrow["{\vp|_{W_k}}", from=2-1, to=2-2]
\end{tikzcd}\]
if we start with $T_i \in V_i$ from the upper left, then since $T_i$ is mapped to $1$ under map $V_i \to W_k$ and $\vp\rst{V_i}(T_i) = 0$ we see $\vp\rst{W_k}$ maps $1$ to $0$.
In the latter case, $\vp\rst{V_i}(F_i) = 0$, so because of the arrow from $V_i$ to $W_k$ we see $\vp\rst{W_k} = 0$ by similar reasons.

Thus, $\vp\rst{W_k} = 0$, so for each $i \in \db{1, M}$, we have $\vp\rst{V_i}(e_i(k)) = 0$ for $1 \le i \le M$ such that $k \in \Omega_i$ because of the arrows from $W_k$ to $V_i$.
However, for each $i \in \db{1, M}$, we have either $\vp\rst{V_i}(T_i) = 0$ or $\vp\rst{V_i}(F_i) = 0$, so from the choice of $e_i(k)$'s we see $\vp\rst{V_i} = 0$.

As a result, $\vp = 0$, i.e., $R_\Theta(x_1, \ldots, x_M)$ is a Schur representation.

Furthermore, we show that $\dim R_\Theta$ is a Schur root no matter whether $\Theta$ can evaluate to true or not.
To see this, just note that we can modify the letters in the clauses of $\Theta$ so that only $X_i$ appear and no $\ob{X_i}$ appear (this can be done by replacing all $\ob{X_i}$ by $X_i$).
Let the modified version of $\Theta$ be $\Theta'$, then $\Theta'$ is able to evaluate to true (by plugging in $X_i = 1$ for all $i$).
However, we have $\dim R_\Theta = \dim R_{\Theta'}$, and the above argument shows that $\dim R_{\Theta'}$ is a Schur root, so $\dim R_\Theta$ is also a Schur root, as desired.

\bigskip

\section{Dimension of $R_{\Theta}(x_1, \ldots, x_M)$ is an imaginary Schur root}
\label{sec: imaginary root}

Let $\Gamma$ be a free abelian group with basis $\br{\al_v}$, labelled by the vertices of the graph $Q$, and we assume that $Q$ has no self loops.
Define a symmetric integer valued bilinear form $(\cdot, \cdot)$ on $\Gamma$ by letting $(\al_v, \al_v) = 2, (\al_v, \al_w) = -\#(\text{edges connecting } v \text{ and } w)$.
Let $r_v$ be the orthogonal reflection of $\Gamma$ with respect to $\al_v$, i.e.,
\[
r_v(\gamma) = \gamma - (\gamma, \al_v) \al_v,
\]
and let $W$ be the group, generated by these reflections.

The element $\dim R = \sum_v (\dim U_v) \al_v \in \Gamma$ is called the \underline{dimension} of the representation $R$ of the quiver $Q$.

Denote by $\Delta_Q$ the set of dimensions of absolutely indecomposable representation of the quiver $Q$ over a finite field $\FF$.
Recall that this set is independent of the orientation of $Q$ and of the field $\FF$, and is called the set of \underline{roots} of the graph $Q$ \cite{K80}.
Recall that the set of roots $\Delta_Q \setminus \br{\al_v}$ is $r_v$-invariant \cite{K80}.
Recall that a root $\al$ is called \underline{imaginary} if and only if its $W$-orbit intersects the fundamental domain $C = \br{\al \in \Gamma \mid (\al, \al_v) \le 0 \text{ for all } v}$.
A root $\al$ is called a \underline{Schur root} if there exists a Schur representation of dimension $\al$.
A Schur root $\al$ is imaginary if and only if $(\al, \al) \le 0$.

In the special case of a quiver $Q_\Theta$ it is natural to take $\Gamma$ a free abelian group with basis
\[
\al_{u_1}, \ldots, \al_{u_M}, \al_{v_1}, \ldots, \al_{v_M}, \al_{w_1}, \ldots, \al_{w_L}.
\]

Then all non-zero inner products between these elements are:
\begin{align*}
    (\al_p, \al_p) = 2, & \text{ for each vertex } p \\
    (\al_{u_i}, \al_{u_{i + 1}}) = -1, & \text{ where } i = 1, \ldots, M \text{ and } u_{M + 1} = u_1; \\
    (\al_{u_i}, \al_{v_i}) = -1, & \text{ where } i = 1, \ldots, M; \\
    (\al_{v_i}, \al_{w_k}) = -2, & \text{ where } i = 1, \ldots, M, k \in \Omega_i.
\end{align*}

We have
\[
\al := \dim R_\Theta(x_1, \ldots, x_M) = \sum_{i = 1}^M \al_{u_i} + \sum_{i = 1}^M n_i \al_{v_i} + \sum_{k = 1}^L \al_{w_k}.
\]

Note that $\al$ is always a Schur root since it admits a Schur representation.
To see this, just note that in the 3-cnf formula, if we replace all $\ob{X_i}$ by $X_i$, then the resulting 3-cnf formula has a solution (let each $X_i$ to be true), and the previous section tells us that the root corresponding to this 3-cnf formula is a Schur root.
However, these two roots are exactly the same, so we see $\al$ is indeed a Schur root.

Now
\[
(\al, \al_{u_i}) = -n_i, (\al, \al_{v_i}) = 1, (\al, \al_{w_k}) = 2 - 2 \sum_{\substack{1 \le i \le M \\ \mathstrut k \in \Omega_i}} n_i = -4 - 2\sum_{\substack{1 \le i \le M \\ \mathstrut k \in \Omega_i}} \abs{\Omega_i}.
\]

Therefore
\[
(\al, \al) = -\sum_{i = 1}^M n_i + \sum_{i = 1}^M n_i - 4 L - 2 \sum_{i = 1}^M \abs{\Omega_i}^2 = - 4 L - 2 \sum_{i = 1}^M \abs{\Omega_i}^2 \le 0,
\]
from which we see $\al$ is an imaginary Schur root.

% Letting $\beta = \prod_{i = 1}^M r_{\al_{v_i}}(\al)$, we see that
% \[
% (\beta, \al_{u_i}) = -n_i, (\beta, \al_{v_i}) = -1, (\beta, \al_{w_k}) = 2 - 2 \sum_{\substack{1 \le i \le M \\ \mathstrut k \in \Omega_i}} n_i,
% \]
% hence $\beta \in C$ and $\al$ is an imaginary root.
% By Section~\ref{sec: eval 1 is schur}, it is an imaginary Schur root.

\section{Construction for an arbitrary finite field $\FF$}

We can generalize the above construction to arbitrary finite field $\FF$.

Here is an example illustrating the construction for $(X_1 \vee X_2 \vee \ob{X_3}) \wedge (X_2 \vee X_3 \vee X_5) \wedge (X_3 \vee \ob{X_4} \vee \ob{X_5})$ (so $M = 5, L = 3$).
The quiver here is just taking the $w$ part of the previous quiver $B := (\abs{\FF} - 1)^3$ times, but we just draw $2$ copies of each $w$ (instead of $B$) to avoid making the diagram too complicated.

All $U_i$'s and $W_k\up{\ell} (1 \le \ell \le B)$'s will have dimension 1, and $\dim V_i = n_i$, where $n_i = B \abs{\Omega_i} + 1$.
All arrows are the same as in the previous construction, except that all $w_k\up{\ell}$ will have arrows coming in and out the same as the previous $w_k$:

% https://q.uiver.app/#q=WzAsMTYsWzMsMCwid18yXnsoMSl9Il0sWzAsMiwid18xXnsoMSl9Il0sWzcsMSwid18zXnsoMSl9Il0sWzIsNCwidl8xIl0sWzMsMywidl8yIl0sWzUsMywidl8zIl0sWzYsNCwidl80Il0sWzQsNSwidl81Il0sWzIsOCwidV8xIl0sWzMsNywidV8yIl0sWzUsNywidV8zIl0sWzYsOCwidV80Il0sWzQsOSwidV81Il0sWzEsMSwid18xXnsoMil9Il0sWzUsMCwid18yXnsoMil9Il0sWzgsMiwid18zXnsoMil9Il0sWzgsOV0sWzksMTBdLFsxMCwxMV0sWzExLDEyXSxbMTIsOF0sWzgsM10sWzksNF0sWzEwLDVdLFsxMSw2XSxbMywxLCIiLDAseyJvZmZzZXQiOi0xfV0sWzEsM10sWzQsMV0sWzEsNCwiIiwwLHsib2Zmc2V0IjotMX1dLFs1LDFdLFsxLDUsIiIsMCx7Im9mZnNldCI6LTF9XSxbNCwwXSxbMCw0LCIiLDAseyJvZmZzZXQiOi0xfV0sWzUsMF0sWzAsNSwiIiwwLHsib2Zmc2V0IjotMX1dLFs3LDBdLFswLDcsIiIsMCx7Im9mZnNldCI6LTF9XSxbNSwyXSxbMiw1LCIiLDAseyJvZmZzZXQiOi0xfV0sWzYsMl0sWzIsNiwiIiwwLHsib2Zmc2V0IjotMX1dLFs3LDJdLFsyLDcsIiIsMCx7Im9mZnNldCI6LTF9XSxbMTIsN10sWzMsMTNdLFsxMywzLCIiLDAseyJvZmZzZXQiOi0xfV0sWzQsMTNdLFsxMyw0LCIiLDAseyJvZmZzZXQiOi0xfV0sWzUsMTNdLFsxMyw1LCIiLDAseyJvZmZzZXQiOi0xfV0sWzQsMTRdLFsxNCw0LCIiLDAseyJvZmZzZXQiOi0xfV0sWzUsMTRdLFsxNCw1LCIiLDAseyJvZmZzZXQiOi0xfV0sWzcsMTRdLFsxNCw3LCIiLDAseyJvZmZzZXQiOi0xfV0sWzUsMTVdLFsxNSw1LCIiLDAseyJvZmZzZXQiOjF9XSxbNywxNV0sWzE1LDcsIiIsMCx7Im9mZnNldCI6LTF9XSxbNiwxNV0sWzE1LDYsIiIsMCx7Im9mZnNldCI6LTF9XV0=
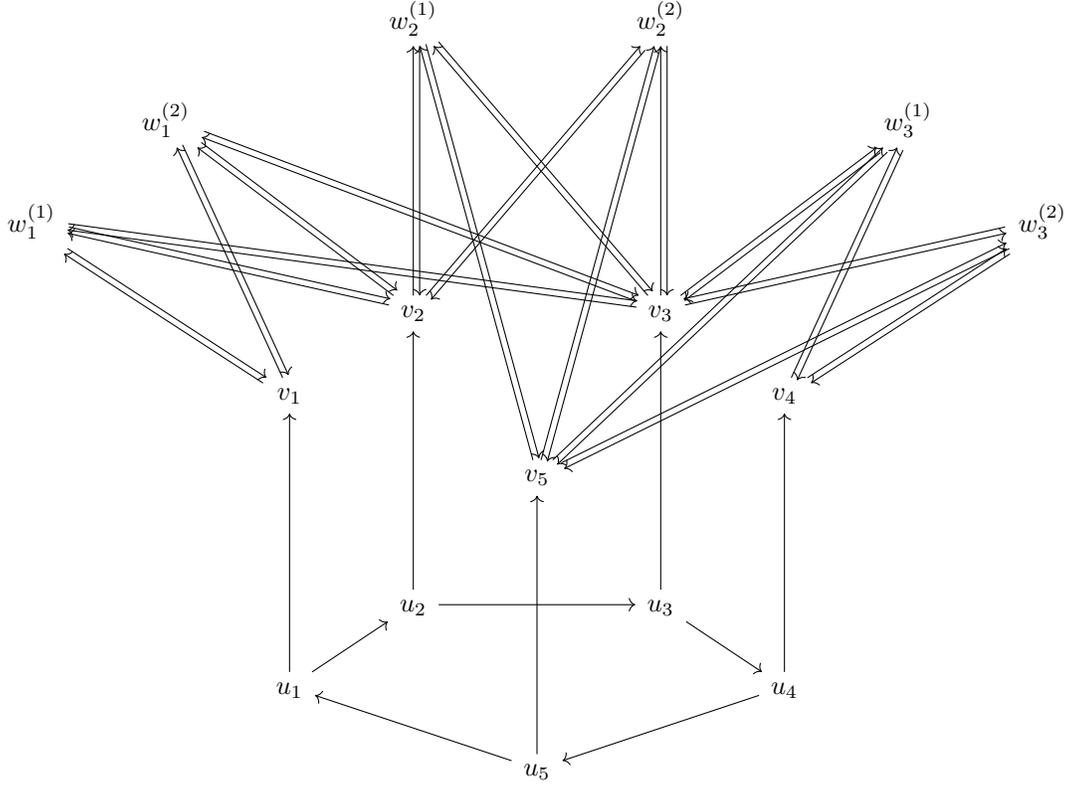
\[\begin{tikzcd}
	&&& {w_2^{(1)}} && {w_2^{(2)}} \\
	& {w_1^{(2)}} &&&&&& {w_3^{(1)}} \\
	{w_1^{(1)}} &&&&&&&& {w_3^{(2)}} \\
	&&& {v_2} && {v_3} \\
	&& {v_1} &&&& {v_4} \\
	&&&& {v_5} \\
	\\
	&&& {u_2} && {u_3} \\
	&& {u_1} &&&& {u_4} \\
	&&&& {u_5}
	\arrow[shift left, from=1-4, to=4-4]
	\arrow[shift left, from=1-4, to=4-6]
	\arrow[shift left, from=1-4, to=6-5]
	\arrow[shift left, from=1-6, to=4-4]
	\arrow[shift left, from=1-6, to=4-6]
	\arrow[shift left, from=1-6, to=6-5]
	\arrow[shift left, from=2-2, to=4-4]
	\arrow[shift left, from=2-2, to=4-6]
	\arrow[shift left, from=2-2, to=5-3]
	\arrow[shift left, from=2-8, to=4-6]
	\arrow[shift left, from=2-8, to=5-7]
	\arrow[shift left, from=2-8, to=6-5]
	\arrow[shift left, from=3-1, to=4-4]
	\arrow[shift left, from=3-1, to=4-6]
	\arrow[from=3-1, to=5-3]
	\arrow[shift right, from=3-9, to=4-6]
	\arrow[shift left, from=3-9, to=5-7]
	\arrow[shift left, from=3-9, to=6-5]
	\arrow[from=4-4, to=1-4]
	\arrow[from=4-4, to=1-6]
	\arrow[from=4-4, to=2-2]
	\arrow[from=4-4, to=3-1]
	\arrow[from=4-6, to=1-4]
	\arrow[from=4-6, to=1-6]
	\arrow[from=4-6, to=2-2]
	\arrow[from=4-6, to=2-8]
	\arrow[from=4-6, to=3-1]
	\arrow[from=4-6, to=3-9]
	\arrow[from=5-3, to=2-2]
	\arrow[shift left, from=5-3, to=3-1]
	\arrow[from=5-7, to=2-8]
	\arrow[from=5-7, to=3-9]
	\arrow[from=6-5, to=1-4]
	\arrow[from=6-5, to=1-6]
	\arrow[from=6-5, to=2-8]
	\arrow[from=6-5, to=3-9]
	\arrow[from=8-4, to=4-4]
	\arrow[from=8-4, to=8-6]
	\arrow[from=8-6, to=4-6]
	\arrow[from=8-6, to=9-7]
	\arrow[from=9-3, to=5-3]
	\arrow[from=9-3, to=8-4]
	\arrow[from=9-7, to=5-7]
	\arrow[from=9-7, to=10-5]
	\arrow[from=10-5, to=6-5]
	\arrow[from=10-5, to=9-3]
\end{tikzcd}\]

We now specify the maps.
The maps $U_i \to V_i$ are still the same, sending $1$ to $(1, X_i, 0, \ldots, 0)$.
The maps $W_k\up{\ell} \to V_i (k \in \Omega_i)$ send $1$ to $e_i(k, \ell) \in V_i$ in such a way that:
\begin{enum}
    \item $\br{e_i(k, \ell) : k \in \Omega_i, 1 \le \ell \le B} \cup \br{(1, c, 0, \ldots, 0)}$ is a basis for $V_i$ for all $c \in \FF$.
    \item $\br{(1, 0, \ldots, 0), (1, c, 0, \ldots, 0)} \cup (\br{e_i(k, \ell) : k \in \Omega_i, 1 \le \ell \le B} - \br{e_i(k', \ell')})$ is a basis for $V_i$ for all $c \in \FF, k' \in \Omega_i, \ell' \in \db{1, \abs{\FF} - 1}$.
\end{enum}
This can be done by choosing one of $e_i(k_0, \ell_0)$ to be $(0, 1, 1, \ldots, 1)$, and the others to be
\[
(0, 0, 1, 0, \ldots, 0), (0, 0, 0, 1, 0, \ldots, 0), \ldots, (0, 0, \ldots, 0, 1).
\]

If $k \in \Omega_i^+$, then we choose the map $V_i \to W_k\up{\ell}$ to send  $(1, 0, \ldots, 0)$ to 0, $(1, 1, 0, \ldots, 0)$ to 1, and all $\br{e_i(k', \ell') : (k', \ell') \ne (k, \ell)}$ to 0.
If $k \in \Omega_i^-$, then the $B$ spaces $W_k\up{\ell}$ admit a one to one correspondence to the triples $(b_1, b_2, b_3) \in (\FF - \br{0})^3$, so if $E_k = \beta_1 \vee \beta_2 \vee \beta_3$ and $\beta_u = \ob{X_i}$ ($1 \le u \le 3$), we choose the map $V_i \to W_k\up{\ell}$ to send $(1, b_u, 0, \ldots, 0)$ to 0, $(1, 0, \ldots, 0)$ to 1, and $\br{e_i(k', \ell') : (k', \ell') \ne (k, \ell)}$ to 0.

We now argue that when $(x_1, \ldots, x_M) \in \FF_2^M$ makes $\Theta$ evaluate to True, then $R = R_\Theta(x_1, \ldots, x_M)$ is a Schur representation, where here $x_1, \ldots, x_M$ are viewed as elements in $\FF$ by the map $\FF_2 \to \FF, 0 \mapsto 0, 1 \mapsto 1$.
To see this, pick any $\vp \in \End R$.
Since $\FF \cdot \id \sub \End R$, after subtracting a suitable multiple of $\id$ we can assume that $\vp\rst{U_1} = 0$.
Then $\vp\rst{U_2} = 0$ by similar reasons as before, and repeating we obtain $\vp\rst{U_i} = 0$ for all $i$.

Because of the arrows from $U_i$ to $V_i$, for each $i \in \db{1, M}$ such that $x_i = 1$, we have $\vp\rst{V_i}(1, 1, 0, \ldots, 0) = 0$, and for each $i \in \db{1, M}$ such that $x_i = 0$, we have $\vp\rst{V_i}(1, 0, \ldots, 0) = 0$ (by similar reasons as before).

For each $k \in \db{1, L}$, there must be some $X_i$ appearing in $E_k$ with $x_i = 1$ or some $\ob{X_i}$ appearing in $E_k$ with $x_i = 0$ (because $(x_1, \ldots, x_M) \in \FF_2^M$ makes $\Theta$ evaluate to 1, so makes each $E_k$ evaluate to 1).
In the former case, $\vp\rst{V_i}(1, 0, \ldots, 0) = 0$, so because of the arrow from $V_i$ to $W_k\up{\ell}$ we see that $\vp\rst{W_k\up{\ell}} = 0$ (by similar reasons as before).
In the latter case, $\vp\rst{V_i}(1, q, \ldots, 0) = 0$ for some nonzero $q \in \FF$, so because of the arrow from $V_i$ to $W_k\up{\ell}$ we see that $\vp\rst{W_k\up{\ell}} = 0$ (by similar reasons as before).

Thus, $\vp\rst{W_k\up{\ell}} = 0$, so for each $i \in \db{1, M}$, we have $\vp\rst{V_i}(e_i(k, \ell)) = 0$ for $1 \le i \le M, 1 \le \ell \le \abs{\FF} - 1$ such that $k \in \Omega_i$ because of the arrows from $W_k\up{\ell}$ to $V_i$ (by similar reasons as before).
However, for each $i \in \db{1, M}$, we have either $\vp\rst{V_i}(1, 1, 0, \ldots, 0) = 0$ or $\vp\rst{V_i}(1, 0, \ldots, 0) = 0$, so from the choice of $e_i(k, \ell)$'s we see $\vp\rst{V_i} = 0$.

As a result, $\vp = 0$, i.e., $R_\Theta(x_1, \ldots, x_M)$ is a Schur representation.

Furthermore, we show that $\dim R_\Theta$ is a Schur root no matter whether $\Theta$ can evaluate to true or not.
To see this, just note that we can modify the letters in the clauses of $\Theta$ so that only $X_i$ appear and no $\ob{X_i}$ appear (this can be done by replacing all $\ob{X_i}$ by $X_i$).
Let the modified version of $\Theta$ be $\Theta'$, then $\Theta'$ is able to evaluate to true (by plugging in $X_i = 1$ for all $i$).
However, we have $\dim R_\Theta = \dim R_{\Theta'}$, and the above argument shows that $\dim R_{\Theta'}$ is a Schur root, so $\dim R_\Theta$ is also a Schur root, as desired.

When $\Theta$ always evaluate to False, then $R_\Theta(x_1, \ldots, x_M)$ will always be decomposable.
To see this, pick any $x_1, \ldots, x_M \in \FF^M$, then there exists $E_k = \beta_{i_1} \vee \beta_{i_2} \vee \beta_{i_3}$ such that
\begin{enum}
    \item Whenever $\beta_{i_j} = \ob{X_{i_j}}$, we have $x_i \ne 0$.
    \item Whenever $\beta_{i_j} = X_{i_j}$, we have $x_i = 0$.
\end{enum}

We pick $\ell$ such that $W_k\up{\ell}$ corresponds to $(y_1, y_2, y_3)$, where $y_j = \begin{cases}
    x_{i_j} & x_{i_j} \ne 0 \\
    1 & x_{i_j} = 0
\end{cases}$.
Then $R_\Theta(x_1, \ldots, x_M) = R_1 \oplus R_2$, where
\[
R_1 = W_k\up{\ell} \oplus \bigoplus_{\substack{i \in \db{1, M} \\ k \in \Omega_i}} \FF e_i(k, \ell) \sub W_k\up{\ell} \oplus \bigoplus_{\substack{i \in \db{1, M} \\ k \in \Omega_i}} V_i.
\]
and
\[
R_2 = \bigoplus_{i = 1}^M U_i \oplus \bigoplus_{i = 1}^M V_i' \oplus \bigoplus_{(k', \ell') \ne (k, \ell)} W_{k'}\up{\ell'} \sub \bigoplus_{i = 1}^M U_i \oplus \bigoplus_{i = 1}^M V_i \oplus \bigoplus_{(k', \ell') \ne (k, \ell)} W_{k'}\up{\ell'},
\]
where
\[
V_i' = \begin{cases}
    V_i & k \not\in \Omega_i \\
    \FF (1, x_i, 0, \ldots, 0) \oplus \Span\br{e_i(k', \ell') : k' \in \Omega_i, 1 \le \ell' \le B, (k', \ell') \ne (k, \ell)} & k \in \Omega_i
\end{cases}.
\]
The reason why this is a valid decomposition is the same as before.

Now if we denote the root by $\al$, then
\begin{align*}
    (\al, \al) & = 2\pr{M + \sum_i n_i^2 + B L} - 2\pr{M + \sum_i n_i + 2 \sum_i n_i B \abs{\Omega_i}} \\
        & = 2\pr{B L - \sum_i n_i B \abs{\Omega_i}} \\
        & = 2B\pr{L - \sum_i \abs{\Omega_i}(B\abs{\Omega_i} + 1)} \\
        & \le 0,
\end{align*}
so $\al$ is an imaginary Schur root.

\section{Acknowledgment}

We are grateful to Mike Sipser for valuable discussions.

% \section{$\al$ is an Imaginary Root}

% \blank Let $\al$ be the root we constructed (depending on $\Theta$).
% Let $e_{u_i}, e_{v_i}, e_{w_k}$ be the generators, then
% \begin{align*}
%     (\al, e_{u_i}) & = -\frac{n_i + 1}{2} \\
%     (\al, e_{v_i}) & = \frac{1}{2} \\
%     (\al, e_{w_k}) & = 1 - \sum_{\substack{i \\ k \in \Omega_i}} (n_i + 1) = -2 - \sum_{\substack{i \\ k \in \Omega_i}} n_i,
% \end{align*}
% and if $\beta$ is the product of reflections of $\al$ through all $v_i$'s (in the Weyl group) (note that the order does not matter), then
% \begin{align*}
%     (\beta, e_{u_i}) & = -\frac{n_i}{2} \\
%     (\beta, e_{v_i}) & = -\frac{1}{2} \\
%     (\beta, e_{w_k}) & = 1 - \sum_{\substack{i \\ k \in \Omega_i}} n_i.
% \end{align*}

% Thus, we see $\al$ is an imaginary root but does not lie in the fundamental chamber.

\end{document}